\documentclass[10pt,a4paper]{article}

\usepackage{amsmath,amsthm}
\usepackage[numbers]{natbib}
\usepackage{microtype}
\usepackage[bitstream-charter]{mathdesign}
\usepackage{bm}
\usepackage[final]{showkeys}
\usepackage{color}

\usepackage[hmargin=3cm,vmargin={3.5cm,4cm}]{geometry}
\allowdisplaybreaks

\title{Discussion on the paper \emph{On Simulation and Properties of the Stable Law} by L.~Devroye and L.~James}

\author{$\text{Mirko D'Ovidio}_1$, $\text{Federico Polito}_2$ \\
	\footnotesize (1) -- Dipartimento di Scienze di Base e Applicate per l'Ingegneria,
		\emph{Sapienza} Universit\`{a} di Roma\\
	\footnotesize Via A. Scarpa 16, 00161 Roma, Italy\\
	\footnotesize Email address: mirko.dovidio@uniroma1.it\\
	\footnotesize (2) -- Dipartimento di Matematica \emph{G.~Peano}, Universit\`{a} degli Studi di Torino\\
	\footnotesize Via Carlo Alberto 10, 10123 Torino, Italy\\
	\footnotesize Email address: federico.polito@unito.it\\
	}

	\date{}

\begin{document}

	\maketitle
	
	\begin{abstract}

		We congratulate the authors for the interesting paper. The reading has been really pleasant and instructive. We discuss briefly only some
		of the interesting results given in \cite{dev} with particular attention to evolution problems. The contribution of the results collected
		in the paper is useful in a more wide class of applications in many areas of applied mathematics. 

	\end{abstract}

		The paper under discussion is a very well-written and interesting review article by Professors Devroye and James
		\citep{dev} dealing with known and lesser-known properties of stable laws,
		with methods of simulation for stable random variates, and with related random variables such as the Mittag--Leffler, Linnik,
		and Lamperti random variables.
		The main aim of the paper is to review and to collect in a single place simple procedures (one-liners) to generates random deviates
		from distributions that are in some way related to stable random variables. This is undoubtedly a very important topic
		at the basis of many techniques in different scientific fields. One can think for example at simulation of stochastic processes,
		generation of pseudo random numbers, cryptography, Monte Carlo and MCMC techniques, and so forth.		 
		A simple example
		in which generation of random deviates that are functions of stable random variables is
		needed, and which can make evident the important of the topic,
		regards the simulation of trajectories of time-fractional point
		processes such as the fractional Poisson process \citep{ors,laskin,mainscalas} or the
		fractional Yule process (fractional pure birth process) \citep{polito} (or in general of renewal processes
		with inter-arrival times distribution related to the stable law). The fractional Poisson process for example
		can be indeed constructed by exploiting its
		renewal structure. Let us thus consider a sequence of iid positive-Linnik distributed random variables $(T_j)_{j=1}^n$ with
		parameter $\mu>0$ (also known
		as Mittag--Leffler random variables in part of the literature, see for example \citet{gerd,pillai,jose})
		which we consider as random inter-arrival times between occurrences of point events. We have that for each $j=1,\dots,n$,
		\begin{align}
			\label{go}
			\mathbb{E} \, e^{sT_j} = \mu/(s^\nu+\mu), \qquad \mathbb{P} \{ T_j \in \mathrm dt \}/\mathrm dt = \mu t^{\nu-1} E_{\nu,\nu}
			(-\mu t^\nu), \qquad \nu \in (0,1], \: \mu>0, \: t > 0.
		\end{align}
		The state probability, that is the probability of attaining level $k$ at time $t$ for the fractional Poisson process $N^\nu(t)$, $t > 0$,
		easily follows from \eqref{go} and reads
		\begin{align}
			\mathbb{P}\{ N^\nu(t) = k \} = (\mu t^\nu)^k E_{\nu,\nu k+1}^{k+1}(-\mu t^\nu), \qquad t > 0, \: k \ge 0,
		\end{align}
		where $E_{\xi,\mu}^\gamma(z)$, for $z,\xi,\mu,\gamma \in \mathbb{C}$, $\Re(\xi)>0$,
		is the three-parameter generalized Mittag--Leffler function \citep{kil}.
		Clearly, the simulation of a trajectory corresponds to the generation of a sequence of independent random variates from positive-Linnik
		distributions. As suggested by the authors (see \citet{dev} but also \citet{dev2,dev3}) a positive Linnik random variate
		can be generated as
		\begin{align}
			\mathcal{E}^{1/\nu} S_\nu,
		\end{align}
		where $\mathcal{E}$ is an exponential random variable of parameter $\mu$ and $S_\nu$ is a completely positively skewed
		stable random variable independent of $\mathcal{E}$. The stable random variable $S_\nu$ can be in turn generated
		by using the classical Kanter algorithm. A simulation of trajectories of the fractional Poisson process using the
		above representation is in fact implemented in \citet{ca}, Section 3.
		
		The paper under discussion
		thoroughly describes simple and useful distributional representations for many different types of random variables, such as
		stable and strictly stable (symmetric and skewed), shifted Cauchy, Lamperti, Linnik and generalized Linnik, and other less known
		representation for related random variables.
		In the following, in order to highlight the usefulness of the relations connecting stable random variables to other related random variables
		we will outline a possible construction of a subordinated Brownian motion time-changed with a specific skewed stable process.
		In fact, for us it seems of particular importance the relation appearing in the paper under discussion \citep{dev} in the section
		entitled ``The strictly stable law: $\alpha>1$''. The authors recall that
		\begin{align}
			\label{first}
			(S_{\alpha,\rho})_+ \overset{\mathcal{L}}{=} (S_{1/\alpha,\alpha\rho})^{-1/\alpha}_+, \qquad \alpha \rho \le 1, \: \alpha (1-\rho) \le 1,
			\: \alpha \in (1,2], \: \rho \in [0,1],
		\end{align}
		where $(S_{\alpha,\rho})_+ = \max(S_{\alpha,\rho},0)$ is the positive part of the strictly stable random variable $S_{\alpha,\rho}$.
		An interesting specific case of which we will make use in the following is when $\rho=1/\alpha$. In this case, \eqref{first}
		reduces to
		\begin{align}
			(S_{\alpha,1/\alpha})_+ \overset{\mathcal{L}}{=} (S_{1/\alpha,1})_+^{-1/\alpha} \overset{\mathcal{L}}{=} S_{1/\alpha}^{-1/\alpha},
			\qquad \alpha \in (1,2]. \label{sssssS}
		\end{align}
				 
		Let $B(t)$, $t>0$, be a Brownian motion with generator $\Delta$. Let $L^\alpha_t$, $t>0$ with $\alpha \in (0,1)$ be the inverse of
		the stable subordinator $S_\alpha(t)$, $t > 0$, independent of $B(t)$. It can be proved that the time-changed
		process $B(L^\alpha_t)$ is the stochastic solution to the fractional equation
		\begin{equation}
			\frac{\partial^\alpha u}{\partial t^\alpha} = \Delta u \quad \text{ in } \quad D \subseteq \mathbb{R} \label{fracPDE},
		\end{equation}
		subject to the Delta initial datum $u_0=\delta$ (see for example \cite{BM2001, nane12}). This is to say that
		\begin{equation}
			\mathbb{P}_x (B(L^\alpha_t)  \in \Lambda) = \int_\Lambda u(x,y,t) \, \mathrm dy, \qquad x \in D,\: t>0,
		\end{equation}
        for some Borel set $\Lambda$. The fractional derivative appearing in \eqref{fracPDE} must be understood in the sense of Caputo.
        The process $L^\alpha_t$ is an inverse process in the sense that
        \begin{equation*}
        	L^\alpha_t = \inf \{s\geq 0\, :\, S_{\alpha}(s) \notin [0, t] \}, \qquad \alpha \in (0,1),
        \end{equation*}
        where $S_\alpha (s) \overset{\mathcal{L}}{=} s^{1/\alpha} S_{\alpha}$ and $S_{\alpha} = S_{\alpha, 1}$, but also in the sense that
        \begin{equation*}
        L^\alpha_t =  M_{\alpha} / t \overset{\mathcal{L}}{=} S^{-\alpha}_{\alpha} / t, \qquad \alpha \in (0,1),
        \end{equation*}
        where $M_\alpha$ is a Mittag--Leffler random variable (see \cite{dev}). 
          
        Let us now consider $\alpha \in (1,2]$. From the fact that 
        \begin{equation*}
			B\left( \frac{1}{t}\right) \overset{\mathcal{L}}{=} \frac{1}{t} B(t)
        \end{equation*}
        we  can write
        \begin{align*}
			B( M_{1/\alpha}/t) \overset{\mathcal{L}}{=} B( S^{-1/\alpha}_{1/\alpha}/t) \overset{\mathcal{L}}{=}
			B(t S_{1/\alpha}^{1/\alpha}) \times S^{-1/{\alpha}}_{1/\alpha} / t
			\overset{\mathcal{L}}{=} B(t S_{1/\alpha}^{1/\alpha}) \times M_{1/\alpha} /t.
		\end{align*}
		In particular, for $\alpha \in (1,2]$, we have that
		\begin{equation}
			B(L^{1/\alpha}_t) \overset{\mathcal{L}}{=} \frac{1}{T_t} B(T_t)
		\end{equation}
		with $T_t = t \, S_{1/\alpha}^{1/\alpha}  \overset{\mathcal{L}}{=} \sqrt[\alpha]{ S_{1/\alpha}(t)}$,
		where $S_{1/\alpha}(t)$ is a stable subordinator such that
		\begin{equation}
			\mathbb{E} \, e^{-\lambda S_{1/\alpha}(t)} = e^{-t \lambda^{1/\alpha}}.
		\end{equation}
		Moreover, it holds that
		\begin{equation}
			B(L^{1/\alpha}_t) \overset{\mathcal{L}}{=}  B\left( \frac{1}{T_t} \right).
		\end{equation}
		From \eqref{sssssS} we also have that
		\begin{equation}
			B(L^{1/\alpha}_t) \overset{\mathcal{L}}{=} B\left( (S_{\alpha,1/\alpha})_+ \, / t\right)
		\end{equation}
		which is really interesting in our view, being $\alpha \in (1,2]$. The considerations above and the results collected in \cite{dev},
		allow us to consider simulation procedures for the solution to more general fractional equations.
		Indeed, we can avoid to deal with inverses and work directly with stable processes. Let us consider a Markov process $X$
		with infinitesimal generator $A$. It is well-known that time-changes of $X$ given by inverses lead to equations of the form
		\begin{equation}
			f(\partial_t) \, u = A \, u,
		\end{equation}
		where $f$ is a well-specified function related to some Bernstein function $g$ (see for example the work \cite{Toaldo} and the references therein).
		For the negative definite operator $A$, it is also known that a time-change given by subordination
		(that is to consider a subordinator as a random time) leads to a time-changed process driven by the equation
		\begin{equation*}
			\partial_t \, u = - g(-A) \, u.
		\end{equation*}
		In this case, the function $g$ is exactly the Bernstein function associated to the time-change, say $\tau_t$, $t>0$, for which
		\begin{equation*}
			\mathbb{E} \, e^{-\lambda \tau_t} = e^{-t g(\lambda)}.
		\end{equation*}
		The semigroup associated with $X_{\tau_t}$ is then a subordinate semigroup.
		
		In our view, the work by Devroye and James turns out to give useful results in many areas of applied mathematics.
		The results collected in this paper has the potential to not only impact the specific applications mentioned above but also the
		more general scenario of applied sciences.

\end{document}